\documentclass[10pt]{amsart}
\usepackage{amssymb}
\usepackage{amsmath}
\usepackage{hyperref}
\usepackage{color}
\newcommand{\lgra}{\longrightarrow}%
\newcommand{\what}{\widehat}%
\newcommand{\wtilde}{\widetilde}%
\newcommand{\Ra}{\mathcal R}%
\newcommand{\R}{\mathbb R}%
\newcommand{\C}{\mathbb C}%
\newcommand{\Q}{\mathbb Q}%
\newcommand{\PP}{\mathcal P}%
\newcommand{\hc}{\mathrm c}%

\newtheorem{theorem}{Theorem}[section]
\newtheorem{thmspecial}{Theorem}

\newtheorem{lemma}[theorem]{Lemma}
\newtheorem{proposition}[theorem]{Proposition}

\theoremstyle{definition}

\theoremstyle{definition}
\newtheorem{remark}[theorem]{Remark}

\numberwithin{equation}{subsection}
\numberwithin{theorem}{subsection}

\begin{document}
\baselineskip14pt

\author[R. P. Sarkar]{Rudra P. Sarkar}
\address[R. P. Sarkar]{Stat-Math Unit, Indian Statistical
Institute, 203 B. T. Rd., Calcutta 700108, India.} 
\email{rudra@isical.ac.in}

\title[Chaotic Dynamics on DR spaces]{Chaotic Dynamics of the heat semigroup on the Damek-Ricci spaces}
\begin{abstract}
The Damek-Ricci  spaces are solvable Lie groups and  noncompact harmonic manifolds. The rank one Riemannian symmetric spaces of noncompact type sits inside it as a thin subclass. In this note we  establish that for any Damek-Ricci space $S$,  the heat semigroup  generated by certain  perturbation of the Laplace-Beltrami operator is {\em chaotic} on the Lorentz spaces $L^{p,q}(S)$, $2<p<\infty, 1\le q<\infty$ and subspace-chaotic on the weak $L^p$-spaces.  We show that both the amount of perturbation and the range of $p$ are sharp. This generalizes a result in \cite{J-W} which proves that under identical conditions, the heat semigroup mentioned above is {\em subspace-chaotic} on the $L^p$-spaces of the symmetric spaces.
\end{abstract}
\maketitle
\section{Introduction and statements of the results}
This article is inspired  by a recent paper of Ji and Web (\cite{J-W}) in which the authors considered  the heat semigroup $T_t, t\ge 0$,  generated by certain  perturbation (which depends on $p$) of the Laplace-Beltrami operator  of a Riemannian symmetric space $X=G/K$ of noncompact type. The authors in \cite{J-W} have shown that $T_t$ is chaotic on the subspace of $K$-invariant functions of $L^p(G/K)$ and hence is {\em subspace-chaotic} on  $L^p(G/K)$  for $2<p<\infty$. This poses a few clear questions: (1)  Is $T_t$  {\em chaotic} on the full space $L^p(G/K)$ with $p$ in the same range?  (2)  Is the amount of perturbation required sharp? (3) Exactly when or for which function space the chaoticity {\em slips} into the subspace-chaoticity? On the other hand, a study of  \cite{J-W} reveals that the parabolic shape of the spectrum of the Laplacian is crucial for this non-Euclidean phenomenon, while the $p$-dependence of the position of the parabolic region justifies the $p$-dependence of the perturbation. This throws some  vindication  that it might be possible to  extend these results to the non-symmetric generalization of the rank one Riemannian symmetric spaces of noncompact type, namely the Damek-Ricci spaces (which are also known as Harmonic $NA$ or $AN$ groups). We shall address these questions. (See the statements below.) We need some preparation before stating the results.

The Damek-Ricci (DR) spaces are solvable Lie groups as well as harmonic manifolds, but very rarely they are  symmetric spaces. Indeed a general DR space  appears as  a counter example to the Lichnerowicz conjecture   (see \cite{DR, DR2}), citing  that there are noncompact  harmonic manifolds which are not symmetric spaces. However  the  rank one Riemannian symmetric spaces of noncompact type  form a  very thin subclass in the set of DR spaces (see \cite{ADY}).
It is well known that such a   symmetric space $X$ is realized as  a quotient space $G/K$ where $G$ is a connected noncompact semisimple Lie group with finite centre and $K$ is a maximal compact subgroup of $G$. Thus $G$ (as well as $K$) has natural left action on  $X$ and functions on $X$ can be realized as right $K$-invariant functions of $G$. One can thus use the full semisimple machinery and in particular the method of decomposing a function in $K$-types to tackle the questions on the function spaces of $X$. The lack of rotation group in a general DR space  is an important difference, which offers fresh difficulties. We  note in this context that  the concept of radiality in a general DR  space is not connected with the group action, which is in contrast with that of the symmetric spaces where a radial function is simply a  $K$-invariant function.

We need the following definitions to proceed:
A strongly continuous semigroup on a Banach space $B$ is a map $T$ from $[0, \infty)$ to the space of all bounded linear operators from $B$ to $B$,
such that    $T(0)=T_0= I$,   the identity operator on $B$;  for all $t,s \ge 0$,  $T_{t + s} = T_t\, T_s$ and  for all $x_0\in B$, $\|T_t x_0 - x_0 \|\to  0$, as $t\to 0$.
The infinitesimal generator $A$ of a strongly continuous semigroup $T$ is defined by    $A \,x=\lim_{t\downarrow 0} \frac1t\,(T_t - I)\,x$, whenever the limit exists and we write $T_t=e^{tA}$.
\begin{enumerate}
\item[(i)] A semigroup of operators $T_t, t\ge 0$  on a Banach space $B$ is {\em hypercyclic} if there exists a $v\in B$ such that $\{T_t v\mid t\ge 0\}$ is dense in $B$.
\item[(ii)] A point $v\in B$ is {\em periodic} for $T_t$, if there exists a  $t>0$ such that $T_t v=v$.

\item[(iii)] The semigroup  $T_t$ is {\em chaotic} if it is hypercyclic and if its  periodic points make a  dense set in $B$.
\item[(iv)] The semigroup $T_t$ is {\em subspace-chaotic} if there is a closed $T_t$-invariant subspace  $V\neq \{0\}$  of $B$   such that $T_t|_V$ is chaotic on $V$.
\end{enumerate}
We have followed \cite{J-W} where chaos is defined in the sense of Devaney (see \cite{Devaney}). For a comprehensive exposition we refer to \cite{J-W}.
Henceforth     $S$  will denote a  Damek-Ricci (DR) space. When a DR space is a rank one Riemannian symmetric space, we shall denote it by $G/K$.
Let $-\Delta$ be the  Laplace-Beltrami operator  on  $S$.
Throughout this article we shall use the notation $c_p$ for $4\rho^2/pp'$ where $p\ge 1$, $p'=p/(p-1)$ and $\rho=Q/2$, $Q$ being the homogenous dimension of $S$ (see section 3).
If the DR space is  a symmetric space then $\rho$ coincides with the half-sum of positive roots, considered as a scaler. We shall assume that $c_\infty=0$.

Purpose of this note is to establish that for any DR space $S$,  the heat semigroup $T_t=e^{-(\Delta-c)t}$ with $c>c_p$,  is {\em chaotic} on the Lorentz spaces $L^{p,q}(S)$, $2<p<\infty, 1\le q<\infty$ and subspace-chaotic on the weak $L^p$-spaces when $2<p\le \infty$.  We show that both the range of $p$ and the  amount of perturbation $c$   are sharp. This generalizes one of the main results in \cite{J-W} which proves that $T_t$ under identical condition is {\em subspace-chaotic} on the $L^p$-spaces of  the symmetric spaces. We recall that the Lorentz spaces are finer subdivisions of the Lebesgue spaces.
Apart from the chaoticity to non-chaoticity, the use of Lorentz spaces  locates another point of degeneracy, where chaoticity changes to subspace-chaoticity. This conforms the paradigm  that the subspace-chaoticity is more stable than chaoticity, as  the results  assert that $T_t$ is at least subspace-chaotic on $L^{p, q}(S)$ for any $1\le q\le \infty$ if and only if  $p>2$ and $c>c_p$.
Our main results are the following. (See section 2 for any unexplained notation.)
\begin{thmspecial} \label{result-S-1}
For $t\ge 0$ and $c\in \R$, let $T_t=e^{-t(\Delta-c)}$. Then,
\begin{enumerate}
\item[(i)] for $2<p<\infty$,  $1\le q<\infty$,   $T_t$ is chaotic on $L^{p,q}(S)$ if  and only if $c>c_p$;

\item[(ii)] for $2<p \le \infty$, $T_t$ is not chaotic on $L^{p, \infty}(S)$ for any $c\in \R$, but subspace-chaotic if  and only if $c>c_p$;

\item[(iii)] for $2<p<\infty$,  $1\le q\le \infty$,  $T_t$ is not hypercyclic and not subspace-chaotic on $L^{p,q}(S)$ and on $L^\infty(S)$ if  $c\le c_p$.
\end{enumerate}
\end{thmspecial}
Note that (ii)  includes the space $L^\infty(S)=L^{\infty, \infty}(S)$.
Part (iii) emphasizes the drastic changes  caused by the amount of perturbation  $c$.
The next  theorem  establishes  the sharpness of the condition $p>2$.
\begin{thmspecial} \label{result-S-2}
For $t\ge 0$ and $c\in \R$, let $T_t=e^{-t(\Delta-c)}$.
\begin{enumerate}
\item[(i)]  For $1<q\le \infty$, $T_t$ is not chaotic or subspace-chaotic on $L^{2, q}(S)$. If $c\le \rho^2$ then
 $T_t$ is not hypercyclic on $L^{2, q}(S)$.

\item[(ii)]  The semigroup $T_t$ is not hypercyclic (hence not chaotic) and not subspace-chaotic on the spaces $L^1(S)$, $L^{2,1}(S)$ and  $L^{p,q}(S)$,
with $1<p<2, 1\le q\le \infty$.
\end{enumerate}
\end{thmspecial}

 Proving the theorems only for the rank one symmetric spaces (which we recall, form a very small subclass of all DR spaces) would be some what simpler as there one can use the compact boundary of the space. Here instead our argument is based on the noncompact boundary.
The paper is organized as follows. The general preliminaries, definitions and results related to chaos are given section 2. In section 3 we give the basic introduction to DR spaces and arrange an array of tools required for the proofs. In section 4 we prove the results stated above. Finally in section 5 we discuss existence of periodic points of the operator $T_t$ in various function spaces.

\section{Notation and Preliminaries}

\subsection{Generalities}
For any $p\in [1, \infty)$, let $p'=p/(p-1)$ and $\gamma_p=(2/p-1)$. The letters $\R$, $\Q$ and $\C$  denote respectively the set of  real numbers, rational numbers and complex numbers.  For
$z\in \C$,  $\Re z$ and $\Im z$  denote respectively the real and
imaginary parts of $z$.   For a set $A$ in a measure space,  $|A|$  denotes the measure of $A$ and for a set $S$ in a topological space, $S^\circ$ denotes its interior.  The letters $C, C_1, C_2$ etc. will be used for positive constants, whose value may change from one line to another.
Occasionally the constants will be suffixed to show their dependencies on important parameters.

\subsection{Lorentz spaces}  We shall  briefly introduce the Lorentz spaces (see \cite{Graf, S-W,
Ray-Sarkar} for details).
Let $(M, m)$ be a $\sigma$-finite nonatomic measure space, $f:M\longrightarrow \C$ be a
measurable function and $p\in [1, \infty)$, $q\in [1, \infty]$. We
define
\begin{equation*}\|f\|^*_{p,q}=\begin{cases}\left(q\int_0^\infty (td_f(t)^{1/p})^q\frac{dt} t\right)^{1/q}
\ \ \ \ \ \  \textup{ if } q<\infty,\\ \\ \sup_{t>0}td_f(t)^{1/p} \ \ \ \ \ \ \ \ \ \ \ \ \ \ \ \ \ \textup{ if } q=\infty,\end{cases}\end{equation*} where
for $\alpha>0$, $d_f(\alpha)=|\{x \mid f(x)>\alpha\}|$ is the distribution function of $f$.
We take
$L^{p,q}(M)$ to be the set of all measurable $f:M\longrightarrow \C$ such that
$\|f\|^*_{p,q}<\infty$. For $1\le p \le \infty$, $L^{p,p}(M)=L^p(M)$
and $\|\cdot\|_{p,p}^\ast=\|\cdot\|_p$.   We note that the Lorentz  ``norm'' $\|\,\cdot\,\|^\ast_{p,q}$ is actually
a quasi-norm  which  makes the space $L^{p,q}(M)$ a quasi Banach
space (see \cite[p. 50]{Graf}). However for $1<p\le \infty$, there
is an equivalent norm $\|\,\cdot\,\|_{p,q}$ through which it is a Banach
space (see \cite[Theorems 3.21, 3.22]{S-W}). We shall slur over this difference, use the notation $\|\cdot\|_{p,q}$ and consider $L^{p,q}(S)$ a Banach space with this norm whenever $p>1$ and we shall not deal with $L^{1, q}$ spaces where $q$ is other than $1$.
The spaces $L^{p,\infty}(M)$ are known as the weak $L^p$-spaces. Thus weak $L^\infty$ space is  same as the $L^\infty$ space.
Following
properties of the Lorentz spaces will be required (see \cite{Graf}). Henceforth for a  Banach space $B$, its dual  space will be denoted by  $B^\ast$.
\begin{enumerate}
\item[(i)] Simple functions are dense in $L^{p,q}(M), 1<p<\infty$, $1\le q<\infty$, but not in $L^{p, \infty}(M)$, $L^\infty(M)$.
\item[(ii)] Unlike $L^{p, q}(M)$ with $q<\infty$, $L^{p, \infty}(M)$ and $L^\infty(M)$ are not separable.
\item[(iii)] If $q_1\le q_2\le \infty$, then $L^{p, q_1}(M)\subset L^{p, q_2}(S)$ and $\|f\|_{p, q_2}\le \|f\|_{p, q_1}$. If $q_2<\infty$ then      $L^{p, q_1}(M)$ is a dense subspace of  $L^{p, q_2}(S)$.
 \item[(iv)] For $1<p, q<\infty$, $(L^{p,q}(M))^\ast=L^{p', q'}(M)$; $(L^{p,1}(M))^\ast=L^{p',\infty}(M)$; $(L^{p, \infty}(M))^\ast=L^{p', 1}(M)\oplus S$ where elements of $S$ are singular functionals (see \cite{Cw}) and   $(L^\infty(M))^\ast=L^{p,1}(S)\oplus M$ where $M$ consists of  certain finitely additive measures.
\end{enumerate}
\subsection{Chaos and hypercyclicity} In section 1 we have defined chaos and hypercyclicity. For a   detailed account we refer to \cite{J-W}. (See also the references therein, in particular \cite{BM, DSW, del-Emam}). Here we shall limit ourselves to what is needed to make the article self-contained. Let $B$ be a  Banach space and $B^\ast$ be its dual space. For a linear operator $A$ on $B$, let $\sigma_{pt}(A)$ be its point spectrum.  For a strongly continuous semigroup of operators (see section 1) $T_t, t\ge 0$ acting on  $B$, let:
\begin{enumerate}
\item[(a)] $B_0=\{x\in B\mid \lim_{t\to \infty} T_t x=0\}$;
\item[(b)] $B_\infty$  be the  set of $x\in B$ such that for each $\varepsilon>0$ there exists $w\in B$ and $t>0$ with $\|w\|<\varepsilon$ and $\|T_tw-x\|<\varepsilon$;
\item[(c)] $B_{\mathrm{Per}}$  be the set of all periodic points in $B$.
\end{enumerate}
Following is a key result proved in \cite{DSW}.
\begin{theorem}  \label{DSW0}
Let $T_t$ denote a strongly continuous semigroup of operators on a separable Banach space $B$. If $B_\infty$ and $B_0$ are dense in $B$, then $T_t$ is hypercyclic.
\end{theorem}
We  also have  the following necessary conditions (see \cite{DSW, del-Emam}) for $T_t$ being hypercyclic/chaotic on $B$:
\begin{proposition}\label{Not-chaotic-hyper} Let $T_t, t\ge 0$ be a  semigroup of operators generated by $A$ in a Banach space $B$.
\begin{enumerate}
\item[(i)] If $T_t$ is chaotic on $B$ then the intersection of the point spectrum of $A$ with $i\R$ is infinite.

\item[(ii)] If $T_t$  is hypercyclic on $B$ then for the adjoint operator $A^*$ of $A$ on the dual space $B^\ast$, $\sigma_{pt}(A^\ast)=\emptyset$.

\item[(iii)] If $T_t$  is hypercyclic on $B$ then for any $\phi\in B^\ast$, $\phi\neq 0$ the orbit $\{T_t^\ast\phi\mid t\ge 0\}$ is unbounded.
\end{enumerate}
\end{proposition}
\begin{remark} \label{non-subspace-chaotic}
It is clear that these conditions above can actually detect  when $T_t$ is not even subspace-chaotic on a Banach space $B$. Precisely, if $\sigma_{pt}(A)\cap i\R$ is finite or $\sigma_{pt}(A^\ast)\neq\emptyset$ or there exists nonzero $\phi\in B^\ast$  such that $\{T_t^\ast\phi\mid t\ge 0\}$ is bounded then $T_t$ is not only non-chaotic, it is non-subspace-chaotic on $B$.
\end{remark}

\section{Damek-Ricci Spaces} To make the article self-contained we shall briefly introduce the DR spaces in this section.  Details can be  retrieved from  \cite{ACB, ADY, Ray-Sarkar, KRS}.
Along the way, we shall also   prepare all technical tools required to  prove the main theorems. While most of these are  known to the experts, it may not be available in this form. In particular Lemma \ref{herz} is new.
\subsection{DR Spaces}
Let $\mathfrak n=\mathfrak v\oplus \mathfrak z$ be a $H$-type
Lie algebra where $\mathfrak v$ and $\mathfrak z$ are vector spaces
over $\R$ of dimensions $m$ and $l$ respectively. Indeed
$\mathfrak z$ is the centre of $\mathfrak n$ and $\mathfrak v$ is
its ortho-complement with respect to the inner product of
$\mathfrak n$. Then we  know that $m$ is even.  The group law of $N=\exp \mathfrak n$ is given by
$$(X, Y). (X', Y')=((X+X', Y+Y'+\frac 12[X, X'])\ \ X\in \mathfrak v, Y\in \mathfrak z.$$
We shall identify
$\mathfrak v$, $\mathfrak z$ and $N$ with $\R^m$, $\R^l$ and
$\R^m\times \R^l$ respectively.  The group  $A=\{a_t=e^t \mid t\in \R\}$  acts on $N$ by nonisotropic
dilation: $\delta_{t}(X,Y)=(e^{t/2}X, e^{t}Y)$.
Let $S=NA=\{(X, Y, a_t)\mid (X, Y)\in N, t\in \R\}$ be
the semidirect product of $N$ and $A$ under the action above. The group law of $S$ becomes:
$$(X, Y, a_t) (X', Y', a_s)=(X+a_{t/2} X', Y+ a_{t} Y'+ \frac  {a_{t/2}}2 [X, X'], a_{t+s}).$$ It then follows that
$\delta_{t}(X,Y)=a_tna_{-t}$, where $n=(X, Y)$.
The Lie group
$S$ is  solvable, connected and simply connected  with
Lie algebra $\mathfrak s=\mathfrak v\oplus\mathfrak z\oplus\R$ and is   nonunimodular. The
homogenous dimension of $S$ is $Q=m/2+l$. For convenience we shall also use the notation $\rho=Q/2$. We note that $\rho$ corresponds to the half-sum of positive roots when $S=G/K$, a
rank one symmetric space of noncompact type.
The group $S$ is equipped with the left-invariant Riemannian
metric $d$ induced by
\begin{equation*}\langle(X, Z, \ell), (X', Z', \ell')\rangle =\langle X, X'\rangle+\langle Z, Z'\rangle+\ell\ell'
\end{equation*} on $\mathfrak s$.
The associated left invariant Haar measure $dx$ on $S$   is given
by
\begin{equation} \int_S f(x)dx=\int_{N\times A}f(na_t)e^{-Q
t}dtdn,
\label{measure-NA}
\end{equation} where $dn(X,Y)=dX\,dY$ and $dX, dY, dt$ are Lebesgue
measures on $\mathfrak v$, $\mathfrak z$ and $\R$ respectively. For an element $x=na_t\in S$, we shall use the notation $A(x)=t$.
We denote the Laplace-Beltrami operator associated to this Riemannian structure by $-\Delta$.

The group $S$ can also be  realized as the unit ball
\begin{equation*}B(\mathfrak  s)=\{(X, Z, \ell)\in \mathfrak  s\mid |X|^2+|Z|^2+\ell^2<1\}\end{equation*}
via a Cayley transform $\mathcal C:S\longrightarrow     B(\mathfrak s)$ (see
\cite[p.~646--647]{ADY} for details). For an element $x\in S$, let
$$|x|=d(C(x), 0)=d(x, e)=\log \frac{1+\|\mathcal C(x)\|} {1-\|\mathcal C(x)\|},$$ where $e$ is the identity element of $S$. In particular $d(a_t, e)=|t|$.

A function $f$ on $S$ is called {\em radial} if for all $x,y\in
S$, $f(x)=f(y)$ if $d(x,e)=d(y,e)$.
 For a function space $\mathcal L(S)$ on $S$ we denote its
subspace of radial functions by $\mathcal L(S)^\#$.
For a suitable function $f$ on $S$ its radialization $\Ra f$ is
defined as
\begin{equation}\Ra f(x)=\int_{S_\nu}f(y)d\sigma_\nu(y),
\label{radialization}
\end{equation} where $\nu=|x|$ and $d\sigma_\nu$ is the
surface measure induced by the left invariant Riemannian metric on
the geodesic sphere $S_\nu=\{y\in S\mid d(y, e)=\nu\}$ normalized
by $\int_{S_\nu}d\sigma_\nu(y)=1$. It is clear that $\Ra f$ is a
radial function and if $f$ is radial then $\Ra f=f$.
We recall the following properties of the operator $\Ra$ (see \cite{DR, ACB}):
\begin{enumerate}
\item $\langle \Ra \phi, \psi\rangle=\langle \phi, \Ra\psi\rangle=\langle \Ra\phi, \Ra\psi\rangle$ for all $\phi, \psi\in C^\infty_c(S)$;
\item $\Ra(\Delta f)=\Delta(\Ra f)$.
\end{enumerate}
Since $|\Ra f|\le \Ra |f|$ and by (1) above, $\int_Sf(x)dx=\int_S \Ra f(x)dx$, we have
$\|\Ra f\|_1\le \|f\|_1$.
Interpolating (\cite[p. 197]{S-W}) with the trivial $L^\infty$-boundedness of $\Ra$ we have,
$$\|\Ra f\|_{p, q}\le \|f\|_{p,q}, 1<p<\infty, 1\le q\le \infty.$$

For two measurable functions $f$ and $g$ on $S$ we define their
convolution as (see \cite[p. 51]{Foll-harmonic}):
\[f\ast g(x)=\int_Sf(y)g(y^{-1}x)dy=\int_Sf(xy^{-1})g(y)e^{QA(y)}dy\] where $e^{-QA(\cdot)}$ is the modular function of $S$.
For a measurable function $g$ on $S$ let $g^\ast(x)=g(x^{-1})$.
 If $g$ is radial then $g^\ast=g$ as $d(x,e)=d(x^{-1}, e)$. It is easy to see that for measurable functions $f, g, h$ on $S$, $\langle f\ast g, h\rangle=\langle f, h\ast g^\ast\rangle$ if both sides make sense.

The Poisson kernel ${\mathcal P}:S\times N\lgra \R$ is defined  by
${\mathcal P}(na_t, n_1)=\wp_{a_t}(n_1^{-1}n)$ where
\begin{equation}\wp_{a_t}(n)=\wp_{a_t}(V, Z)=C a_t^{Q}\left(\left(a_t+\frac{|V|^2}4\right)^2+|Z|^2\right)^{-Q},\,\, n=(V,
Z)\in N. \label{poisson}\end{equation}  In particular
\begin{equation} \wp_1(n)=\wp_{a_0}(n)=C [(1+\frac{|V|^2}4)^2 +|Z|^2]^{-Q}. \label{poisson-1}
\end{equation} The  value of $C$ is
adjusted so that $\int_N \wp_1(n)dn=1$ (see
\cite[(2.6)]{ACB}). We note that $\wp_{a_t}(n)=\wp_{a_t}(n^{-1})=\wp_1(a_{-t}na_t)e^{-Qt}=\wp_1(\delta_{-t}(n))e^{-Qt}$. The complex power of Poisson kernel $\mathcal P_\lambda$ is  defined by
\[\PP_\lambda(x, n)=\PP(x, n)^{\frac 12-\frac{i\lambda}Q}.\] Then  for each fixed $n\in N$, $\Delta \PP_\lambda(x,n)=(\lambda^2+\rho^2)\PP_\lambda(x, n)$.
The Poisson transform of a function $F$ on $N$ is defined as (see
\cite{ACB})
\begin{equation*}{\mathfrak P}_{\lambda}F(x)=\int_N F(n)\PP_\lambda (x, n)dn.\end{equation*} It follows that
$\Delta\mathfrak P_\lambda F= (\lambda^2+\rho^2) \mathfrak P_\lambda F$.
For $\lambda\in \C$, the elementary spherical function $\phi_\lambda$ is given by
\[\phi_\lambda(x)=\int_N\PP_{\lambda}(x, n)\PP_{-\lambda}(e, n) dn=\int_N\PP_{-\lambda}(x, n)\PP_{\lambda}(e, n) dn.\]
and we have  (see \cite[Prop 4.2]{ACB}),
\begin{align}\phi_\lambda(x^{-1}y)&=\int_N \mathcal P_\lambda(x, n)\mathcal P_{-\lambda}(y, n) dn.
\label{acb-4.2}
\end{align}
It follows that
$\phi_\lambda$ is a radial eigenfunction of  $\Delta$  with eigenvalue $(\lambda^2+\rho^2)$
satisfying $\phi_\lambda(x)=\phi_{-\lambda}(x),
\phi_\lambda(x)=\phi_\lambda(x^{-1})$, $\phi_{-i\rho}\equiv 1$ and $\phi_\lambda(e)=1$.
We have the following asymptotic estimate of $\phi_\lambda$ (see \cite{ADY}):
\begin{equation}|\phi_{\alpha+i\gamma_p\rho}(x)|\asymp e^{-(2\rho/p')  |x|}, \ \ \alpha\in \R, 0<p< 2.\label{exact-estimate-1}\end{equation}
The estimate above degenerates when $p=2$, i.e. when $\gamma_p=0$ and in this case we have  $\phi_0(x)\asymp (1+|x|) e^{-\rho |x|}$.
If $\lambda\in \R$, $\lambda\neq 0$ and $t\ge 1$ then the Harish-Chandra series for $\phi_\lambda$ implies,
\begin{equation} \phi_\lambda(a_t)=e^{-\rho t}[\hc(\lambda)e^{i\lambda t}+\hc(-\lambda)e^{-i\lambda t}+E(\lambda,
t)], \text{ where } |E(\lambda, t)|\le C_{\lambda} e^{-2t},
 \label{phi-delta-asym}
\end{equation} where $\hc(\lambda)$ is the Harish-Chandra $\hc$-function.
See \cite[(3.11)]{Ion-Pois-1}) for  a proof of the above for the symmetric spaces. The proof works
{\em mutatis mutandis} for general Damek-Ricci spaces.
Let $C_0(S)$ be the space of continuous functions vanishing at infinity with supremum norm. Then $C_0(S)$ is a separable Banach space and a non-dense subspace of $L^\infty(S)$.
For $p\ge 1$, let \[S_p=S_{p'}=\{z\in \C\mid |\Im z|\le |\gamma_p|\rho\}.\] Let
$S_p^\circ$ and $\partial S_p$ respectively be the interior and the
boundary of the strip $S_p$.
We recall that (see \cite{Ray-Sarkar, KRS, KRS-2}):
\begin{enumerate}
\item[(a)] $\phi_\lambda\in L^\infty(S)$ if and only if $\lambda\in S_1$;
\item[(b)] $\phi_\lambda\in C_0(S)$ if and only if $\lambda\in S_1^\circ$;
\item[(c)] for $1<p<2$,  $1\le q<\infty$, $\phi_\lambda\in L^{p', q}(S)$ if and only if $\lambda\in S_p^\circ$;
\item[(d)] for  $1<p<2$,   $\phi_\lambda\in L^{p', \infty}(S)$ if and only if $\lambda\in S_p$;
\item[(e)] for $\lambda\in S_2=\R$, $\phi_\lambda\in L^{2, \infty}(S)$ if and only if $\lambda\neq 0$;
\item[(f)] for $\lambda\in S_2=\R$, $\phi_\lambda\not\in L^{2, q}(S)$ for any $q<\infty$.
\end{enumerate}

The spherical Fourier transform of a function $f$ is defined by $\what{f}(\lambda)=\int_Sf(x) \phi_\lambda(x)dx$ whenever the integral converges. The estimates of $\phi_\lambda$ given above determines the domain of the spherical Fourier transform of functions in different Lebesgue and Lorentz spaces. Indeed,
\begin{enumerate}
\item[(a)] For $f\in L^1(S)$, $\what{f}$ extends as an analytic function on $S_1^\circ$ which is continuous on its boundary; for  a complex measure $\mu$ on $S$, $\mu(\lambda)=\int \phi_\lambda(x)d\mu(x)$ behaves the same way;
\item[(b)] for $f\in L^{p, q}(S)$, with $1<p<2$ and $1<q \le\infty$, $\what{f}$ extends as an analytic function on $S_p^\circ$;
\item[(c)] for $f\in L^{p, 1}(S)$, with $1<p<2$, $\what{f}$ extends as an analytic function on $S_p^\circ$ which is continuous on its boundary;
\item[(d)] for $f\in L^{2,1}(S)$, $\what{f}$ is  continuous  on nonzero real numbers.
\end{enumerate}
For a measurable  function $f$ on $S$ following \cite{ACB} we define its  (Helgason-type) Fourier
transform  by $$\wtilde{f}(\lambda,
n)=\int_Sf(x)\PP_\lambda(x, n)dx,$$ whenever the integral
converges. If $g$ is radial  then,
$\wtilde{f\ast g}(\lambda, n)=\wtilde{f}(\lambda, n)\what{g}(\lambda)$ whenever both sides make sense.

For  $f\in L^{p,q}(S), 1<p<2, 1\le q\le \infty$ (respectively $f\in L^1(S)$), $\lambda\mapsto \wtilde{f}(\lambda, n)$ is a holomorphic function on $S_p^\circ$ (respectively on $S_1^\circ$) for every fixed $n\in N_1$ where $N_1$ is a subset of $N$ of full measure (see \cite[Theorem 3.4, Theorem 5.4]{Ray-Sarkar}). The argument in \cite{Ray-Sarkar} also shows that if $\mu$ is a (bounded) complex measure on $S$, then $\wtilde{\mu}(\lambda, n)$ is a holomorphic function on $S_1^\circ$ for every fixed $n\in N_1$ for $N_1$ as above. The inversion formula is the following (see \cite{Ray-Sarkar}). We recall that $|\hc(\lambda)|^{-2}$ is the Harish-Chandra Plancherel measure.
\begin{proposition}
Let $f\in L^{p, q}(S)$, $p\in (1, 2)$, $q\ge 1$ or $f\in
L^1(S)\cup L^2(S)$. If $\wtilde{f}\in L^1(N\times \R,
|\hc(\lambda)|^{-2}\, d\lambda\, dn)$, then for almost every $x\in
S$
\begin{equation*}f(x)=C\int_{N\times \R}\wtilde{f}(\lambda,n) {\mathcal
P}_{-\lambda} (x,n)|\hc(\lambda)|^{-2}d\lambda\, dn.\end{equation*}
\label{inversion}
\end{proposition}
If $f$ is radial then  $\wtilde{f}(\lambda, n)=\what{f}(\lambda)\PP_\lambda(e, n)$ and  the inversion formula reduces to
\begin{equation}f(x)=C\int_{\R}\what{f}(\lambda) \phi_\lambda(x)|\hc(\lambda)|^{-2}d\lambda.\label{inversion-radial}\end{equation}

We have the following estimate of the Fourier transform vis-a-vis the Poisson transform (see \cite[Theorem 1.1]{KRS}) which  will be used in the proof Theorem \ref{result-S-1}. For $p>2$, $p'<r<p$, $\alpha\in \R$ and $f\in L^{p', \infty}(S)$,
\begin{equation}\|\wtilde{f}(\alpha+i\gamma_r\rho, \cdot)\|_{L^r(N)}\le C \|f\|_{p', \infty}
\label{restriction}
\end{equation}
which by duality is equivalent to (for any $F\in L^{r'}(N)$),
\begin{equation}
 \|\mathfrak P_{\alpha+i\gamma_r\rho}F\|_{p,1}\le C\|F\|_{L^{r'}(N)}.
 \label{extension}
 \end{equation}

We conclude this subsection defining the Harish-Chandra Schwartz spaces on $S$.
For $1\le p\le 2$, the $L^p$-Schwartz space $C^p(S)$ is defined (see
\cite{ADY, di-B}) as the set of $C^\infty$-functions on $S$ such that
$$\gamma_{r,D}(f)=\sup_{x\in S}|Df(x)|
\phi_0^{-2/p}(1+|x|)^r<\infty,$$ for all nonnegative integers $r$
and left invariant differential operators $D$ on $S$.
We recall that (see \cite{KRS-2}) for $1\le p\le 2$, $C^p(S)$ is dense in $L^{p,1}(S)$ and hence in $L^{p,q}(S)$ for $1\le q<\infty$, but not in $L^{p, \infty}(S)$.

\subsection{Herz's criterion}
Let $\mu$ be a nonnegative radial finite measure. We consider the right  convolution operator $T_\mu$ defined on the measurable functions on $S$  by  $T_\mu: f\mapsto f\ast \mu$, whenever it makes sense. We have  the following  {\em Herz's criterion} for the Lorentz spaces.
In  \cite[Theorem 3.3]{ADY}) (see also \cite{Cow-Herz},\cite[Theorem 3.2]{CGM1},  \cite[Proposition 4.1]{KRS}) a more general result is obtained  for the Lebesgue spaces on $S$.

\begin{lemma} \label{herz} Let $1<p<\infty, 1\le q\le \infty$ be fixed. If a nonnegative radial measure $\mu$ satisfies $\int_S \phi_{i\gamma_p\rho}(x)d\mu(x)<\infty$, then $T_\mu$ is a bounded operator from $L^{p,q}(S)$ to itself and the  operator norm of $T_\mu$  satisfies  $\|T_\mu\|_{L^{p,q}\to L^{p,q}}\le \int_S \phi_{i\gamma_p\rho}(x)d\mu(x)$.
\end{lemma}
\begin{proof}
For $x\in S$, we define $(R_p(x)f)(y)=e^{-Q/p A(x)}f(yx)$ where for $x=na_t$, $A(x)=t$. We shall show that $\|R_p(x)f\|_{p,q}=\|f\|_{p,q}$ for any fixed $x\in S$.
Indeed for any $\alpha>0$,
 \begin{align*}
 d_{R_p(x)f}(\alpha)&=|\{y\in S\mid |(R_p(x)f)(y)|>\alpha \}|\\
 &=|\{y\in S\mid |f(yx)|>\alpha e^{Q/p A(x)}\}|\\
 &=|\{yx\in S\mid |f(yx)|>\alpha e^{Q/p A(x)}\}| e^{QA(x)}\\
 &=d_f(\alpha e^{Q/p A(x)})e^{QA(x)}.
 \end{align*}
Therefore if $q<\infty$ then,
 \begin{align*}\|R_p(x)f\|_{p,q}^q&=C\int_0^\infty \alpha^{q-1} d_{R_p(x)f}(\alpha)^{q/p}d\alpha\\
 &=C\int_0^\infty \alpha^{q-1} [d_f(\alpha e^{Q/p A(x)}]^{q/p}\; e^{Q qA(x)/p}\, d\alpha  \\
 &=C\int_0^\infty (\alpha e^{Q/p A(x)})^{q-1} (d_f(\alpha e^{Q/p A(x)}))^{q/p} d(\alpha e^{Q/p A(x)}) \\ 
 &=\|f\|_{p,q}^q
 \end{align*}
 and if $q=\infty$ then,
\begin{align*}\|R_p(x)f\|_{p, \infty}&=\sup_{\alpha>0}\alpha d_{R_p(x)f}(\alpha)^{1/p}\\
&=\sup_{\alpha>0}\alpha e^{Q/p A(x)} d_f(\alpha e^{Q/p A(x)})^{1/p}\\
&=\sup_{\beta>0}\beta\, d_f(\beta)^{1/p}=\|f\|_{p, \infty}.\end{align*}

We note that  $T_\mu(f)(y)=\int_S f(yz^{-1})e^{QA(z)}d\mu(z)$.
For  $f\in L^{p,q}(S), g\in L^{p',q'}(S)$ (taking $q'=1$ when $q=\infty$) we have,
\begin{align*}
\langle T_\mu f, g\rangle &=\int_Sf\ast \mu(y)g(y)dy\\
&=\int_S\int_S f(yz^{-1}) e^{QA(z)} d\mu(z) g(y)dy\\
&=\int_S e^{Q/p' A(z)} \left(\int_S (R_p(z^{-1})f)(y) g(y)dy\right) d\mu(z)\\
&\le \int_S e^{Q/p' A(z)}d\mu(z) \|f\|_{p,q}\|g\|_{p',q'}.
\end{align*} We have used the fact that $A(z^{-1})=-A(z)$ in one of the steps above.
Since $\mu$ is radial  and  $\mathcal R(e^{Q/p' A(\cdot)})=\phi_{i\gamma_p\rho}$ (see \cite[3.11]{ADY}) we have (see \cite[p. 70]{Graf}),
 \[\|T_\mu\|_{L^{p,q}-L^{p,q}}\le \int_S e^{Q/p' A(z)} d\mu(z)=\int_S \phi_{i\gamma_p\rho}(z) d\mu(z)=\what{\mu}(i\gamma_p\rho).\qedhere\]
  \end{proof}
Through similar steps one can also show that if a nonnegative radial measure $\mu$ satisfies $\int_S d\mu(x)<\infty$, then \[\|T_\mu\|_{L^1\to L^1}\le \int_S d\mu(x)\quad \text{and}\quad \|T_\mu\|_{L^\infty\to L^\infty}\le \int_S d\mu(x).\]
\subsection{Spectrum of $\Delta$} We recall that for $p\ge 1$,  $\gamma_p=(2/p-1)$ and $S_p=S_{p'}=\{z\in \C\mid |\Im z|\le |\gamma_p|\rho\}$.
Under the map $\Lambda: z\mapsto z^2+\rho^2$, $S_p$ is mapped to a parabolic region
 in the complex plane which we shall denote by $P_p$. Precisely $P_p=\Lambda(S_p)$. We note that $P_p=P_{p'}$ and if $p=2$ then $S_p=\R$ and
 the Parabolic region reduces to the ray $[\rho^2, \infty)$.  For $p\neq 2$, let $S_p^\circ$ and $P_p^\circ$ be the interiors of $S_p$  and $P_p$ respectively.
 We enlist here some information related to the spectrum of $\Delta$  which will be useful for proving our main results.

Suppose that a nonzero measurable function $u$ satisfies  $\Delta u=(\lambda^2+\rho^2)u$ for some $\lambda\in \C$. We assume that for a point $x_0\in S$, $u(x_0)\neq 0$. Let $\ell_x f$ be the left translation of a function $f$ by $x\in S$.
Since $\Delta$ commutes with the radialization operator and translations, $\Ra (\ell_{x_0}u)$  is a radial eigenfunction with the same eigenvalue $\lambda^2+\rho^2$. Hence
$\Ra (\ell_{x_0}u)=C\phi_\lambda$ (see \cite[2.5]{ADY}).   As $\Ra (\ell_{x_0}u)(e)=\ell_{x_0}u(e)=u(x_0)$ and $\phi_\lambda(e)=1$ we have $\Ra (\ell_{x_0}u)=u(x_0)\phi_\lambda$.
Thus from the $L^{p,q}$-properties of $\phi_\lambda$ given above one can determine  the $L^{p,q}$-point spectrum of $\Delta$. Precisely,
 \begin{enumerate}
 \item[(i)] if $2<p<\infty, 1\le q<\infty$, then
  $P_p^\circ$ (respectively $P_p$) is  the  $L^{p, q}$-point spectrum (respectively the $L^{p, \infty}$-point spectrum);
  \item[(ii)] $P_1$ is the $L^\infty$-point spectrum and $P_1^\circ$ is the $C_0$-point spectrum;
  \item[(iii)] for $1<p<2, 1\le q\le \infty$, the  $L^{p, q}$-point spectrum and the $L^1$-point spectrum  are empty;
\item[(iv)] $(\rho^2, \infty)$ is the $L^{2, \infty}$-point
   spectrum;
  \item[(v)] for $1\le q<\infty$,  the  $L^{2, q}$-point spectrum is empty.
\end{enumerate}

It is known that for $1\le p\le \infty$, the $L^p$-spectrum of $\Delta$ is $P_p$ (see \cite[Cor. 4.18]{ADY}, see also \cite{Ank2, Lo-Ry-2, Tay}).
We restrict  to the range $\{1<p<2\}\cup\{2<p<\infty\}$ and  take a $\lambda\in \C\setminus P_p$. Then we can choose $p_1, p_2$ satisfying $1<p_1<p<p_2<2$ or $2<p_1<p<p_2$ so that $\lambda\not\in P_{p_1}\cup P_{p_2}$. Therefore $\lambda$ is in the $L^{p_1}$-resolvent set as well as in the $L^{p_2}$-resolvent set of $\Delta$. Using interpolation (\cite[p. 197]{S-W}) we conclude that  $\lambda$ is in the $L^{p, q}$-resolvent set of $\Delta$ for any $1\le q\le \infty$. Thus for  $1<p<\infty$ with $p\neq 2$ and $1\le q\le \infty$, the $L^{p, q}$-spectrum of $\Delta$ is a subset of $P_p$. In particular the $L^{p, \infty}$-spectrum of $\Delta$ is $P_p$ which is also the $L^{p, \infty}$-point spectrum   when $p>2$ as mentioned above. We conclude, noting that for an operator $A$ on a Banach space $B$, the spectrum of $T_t=e^{-A t}$ is  not necessarily in  one-to-one correspondence with the spectrum of $A$. (See \cite{DSW}).

\subsection{Heat kernel and the semigroup $T_t$}
Let $h_t$ be the {\em heat kernel} which is defined as a radial function in the Harish-Chandra  Schwartz space $C^p(S)$, $1\le p\le 2$,  by
prescribing its spherical Fourier transform $\what{h_t}(\lambda)=e^{-t(\lambda^2+\rho^2)}$ for all $\lambda\in \C$ (\cite[(5.4), (5.5)]{ADY}).
For any $c\in \R$ and $T_t=e^{-t(\Delta-c)}$,  $T_tf=e^{ct} e^{-t\Delta} f=e^{ct} f\ast h_t$ for any suitable function $f$.
We recall (see \cite[5.50]{ADY}) that the heat maximal operator on $S$, $Mf(x)=\sup_{t>0}|e^{-t\Delta}f(x)|$ is weak type $L^1-L^1$ and strong type $L^\infty-L^\infty$ and hence  by interpolation  (\cite[p. 197]{S-W}), it is strong type $L^{p,q}-L^{p,q}$ for $1<p<\infty$, $1\le q<\infty$. From the fact that $f\ast h_t(x)\to f(x)$  in $C^p(S)$ and that $C^p(S)$ is dense in $L^{p, q}(S)$, the standard method of maximal function yields  $f\ast h_t\to f$ in $L^{p,q}(S)$  as $t\to 0$.  From this it is easy to see that $\|T_tf-f\|_{p,q}\to 0$
 as $t\to 0$ for any $f\in L^{p,q}(S)$.  That is $T_t$ is a strongly continuous semigroup on $L^{p,q}(S)$.
It is easy to verify that  $T_t$ is strongly continuous on  $C_0(S)$.
It is known that $f\ast h_t$ does not converge to $f$ in $L^\infty$ and the same is true for the weak-$L^p$ spaces.

\section{Proof of Theorem A and Theorem B}
\begin{proof}[proof of Theorem \ref{result-S-1}]
(i)  It suffices  to show that $T_t$ is chaotic on $L^{p,1}(S)$, since  for $1\le q<\infty$, $L^{p,1}(S)$ is a dense subspace of $L^{p,q}(S)$. For this proof  $B=L^{p,1}(S)$.

 From the description of $P_p$ (see section 3.3) and the condition $c>c_p$, the following conclusions are immediate: $\Omega_p=(P_p^\circ-c)\cap \{z\mid \Im z>0\}$, $\Omega_p^+=\Omega_p\cap\{z\in \C\mid \Re z>0\}$ and $\Omega_p^-=\Omega_p\cap\{z\in \C\mid \Re z<0\}$ are connected non-empty open sets ; $\Omega_p$ intersects $i\R$ in a nondegenerate line segment. For $z\in \Omega_p$ we define the map  $\Gamma(z)=\sqrt{z-\rho^2+c}$ taking an analytic branch. Then  $\Gamma(z)\in S_p^\circ$ and  $\Gamma:\Omega_p \to S_p^\circ$ is  holomorphic. Hence $\Gamma(\Omega_p^+)$ and $\Gamma(\Omega_p^-)$  are connected and (by the open mapping theorem)  open sets in $S_p$.

Since for $z\in \Omega_p$, $\Gamma(z)\in S_p^\circ$, we have $\Gamma(z)=\alpha+i\gamma_r\rho$ for $\alpha\in \R$ and for some $r$ satisfying $p'<r<p$. Thus every $z\in \Omega_p$ determines an unique $r=r(z)\in (p', p)$, precisely by  $r=r(z)=2\rho(\Im \Gamma(z)+\rho)^{-1}$.

We define, \[U_1=\{\mathfrak P_{\Gamma(z)}F\mid  z\in \Omega_p^+, F\in L^{r'}(N), \text{ where } r=r(z)\}.\] By (\ref{extension}),  $U_1\subset L^{p,1}(S)$, since $r=r(z)$. It is also clear (see section 3) that  elements of $U_1$ satisfy $(\Delta-c)\, \mathfrak P_{\Gamma(z)}\,F=z\, \mathfrak P_{\Gamma(z)}F$ and hence  $T_t\, \mathfrak P_{\Gamma(z)}\, F=e^{-tz}\, \mathfrak P_{\Gamma(z)}F$. The condition $\Re z>0$ ensures that $T_t\psi\to 0$ as $t\to \infty$ for any $\psi\in U_1$. Thus $U_1\subset B_0$. Since $B_0$ is a vector space $\mathrm{span} (U_1)\subset B_0$.

 We assume  that a function $f\in L^{p', \infty}(S)$ annihilates $U_1$. We fix a $z\in \Omega_p^+$  and consider the corresponding  elements  $\mathfrak P_{\Gamma(z)}F$ of $U_1$. Then   by our assumption $\int_S f(x)\mathfrak P_{\Gamma(z)}F(x)dx=0$. This  implies that $\int_N \wtilde{f}(\alpha+i\gamma_r\rho, n)F(n)dn=0$ where $\Gamma(z)=\alpha+i\gamma_r\rho$, $r=r(z)$.  Noting that by (\ref{restriction}),   $n\mapsto \wtilde{f}(\alpha+i\gamma_r\rho, n)$ is in $L^r(N)$, and as  $F$ is an arbitrary function in $L^{r'}(N)$ (where $r=r(z)$), we conclude that  $\wtilde{f}(\alpha+i\gamma_r\rho, \cdot) \equiv 0$. In this way we can show that for any $\lambda\in \Gamma(\Omega_p^+)$, $\wtilde{f}(\lambda, \cdot)\equiv 0$.
  Recalling  that for almost every fixed $n$, $\lambda\mapsto\wtilde{f}(\lambda, n)$ is a holomorphic function on $S_p^\circ$, and that $\Gamma(\Omega_p^+)$ is an open set, we conclude that $\wtilde{f}(\lambda, n)=0$ for all $\lambda\in S_p^\circ$ and  for almost every $n\in N$. By (\ref{inversion}) this implies that $f=0$. This shows that $\mathrm{span} (U_1)$ and hence $B_0$ is dense in $L^{p,1}(S)$.

Next  we define:
  \[ U_2=\{\mathfrak P_{\Gamma(z)}F\mid  z\in \Omega_p^-, F\in L^{r'}(N) \text{ where } r=r(z)\}\text{ and }\]
  \[U_3=\{\mathfrak P_{\Gamma(z)}F\mid z\in \Omega_p\cap i\Q, F\in L^{r'}(N)\text{ where } r=r(z)\}.\]
 Like $U_1$, the sets   $U_2$ and $U_3$ are also  subsets of $L^{p,1}(S)$ and the elements of $U_2$ and $U_3$ satisfy $T_t\mathfrak P_{\Gamma(z)}F=e^{-tz}\mathfrak P_{\Gamma(z)}F$. Let  $\{v_{z_1}, v_{z_2}, \dots, v_{z_n}\}$ be a finite subcollection of elements of $U_2$  with $(\Delta-c)\, v_{z_k}=z_k v_{z_k}$. We take a  $\C$-linear combination of these elements of $U_2$:  \[g=\sum_{k=1}^n a_k v_{z_k}=T_t\left(\sum_{k=1}^n a_k\, e^{z_k t}\, v_{z_k}\right).\] Since $\Re z_k<0$, for any $\varepsilon>0$ we have a  suitable $t>0$, $w=\sum_{k=1}^n a_ke^{z_k t}v_{z_k}$ satisfies $\|w\|_{p,1}<\varepsilon$ and it is clear that $\|T_tw-g\|_{p,1}<\varepsilon$. This shows that   $\mathrm{span}(U_2)\subset B_\infty$. Once we note that $\Gamma(\Omega_p^-)$ is open, the denseness of $U_2$ follows through  the same  argument used for showing denseness of $U_1$. Thus $B_\infty$ is dense.

Finally  it is clear that $\mathrm{span}(U_3)\subset B_{\mathrm{Per}}$. Noting that the set $\Omega_p\cap i\Q$ has limit points we can apply again  almost a similar argument applied for $U_1$ and $U_2$ to show that  $\mathrm{span}(U_3)$ is dense. This proves the  assertion.

To prove the converse, we notice that if $c\le c_p$ then $(P_p-c)\cap i\R$ has at most one point. Therefore by Proposition~\ref{Not-chaotic-hyper}~(i),  $T_t$ is not chaotic in this case.
\vspace{.1in}

(ii) Since $T_t$ is not strongly continuous on $L^\infty(S)$ and on $L^{p, \infty}(S)$, it is easy to see that $T_t$ cannot  be hypercyclic and hence not chaotic on these spaces.
   Since  $c>c_p$, by (i) above,  $T_t$ is chaotic on $L^{p}(S)$ which is a subspace of $L^{p, \infty}(S)$. Therefore $T_t$ is subspace-chaotic on $L^{p, \infty}(S)$ when $c>c_p$.

We shall show  that if $c>c_\infty=0$, then $T_t$ is chaotic on $C_0(S)$ which is a subspace of $L^\infty(S)$. For this proof $B=C_0(S)$. We fix a  $c>0$ and define the sets, $\Omega_\infty$, $\Omega_{\infty}^+$ and $\Omega_{\infty}^-$ putting $p=\infty$ in the definition of $\Omega_p$ given in (i). Then $\Omega_\infty$ is a connected open set which  intersects $i\R$ in a nondegenerate  line segment and $\Gamma(\Omega_\infty)\subset S_1^\circ$, where the  function $\Gamma$ is as defined in (i).
It is clear that  for $z\in \Omega_\infty$, $\psi(z)=\phi_{\Gamma(z)}$ (the elementary spherical function $\phi_\lambda, \lambda=\Gamma(z)$) is an eigenfunction of $\Delta-c$ with eigenvalue $z$ and $\psi(z)\in C_0(S)$ (see section 3).

Let $V_1=\{\ell_y(\psi(z))(x)\mid z\in \Omega_\infty^+,  y\in S\}$, where  $\ell_y$ is the left translation;  $\ell_y (\psi(z))(x)=\phi_{\sqrt{z-\rho^2+c}}(y^{-1}x)$.
Then $\ell_y (\psi(z))(x)$ is also an eigenfunction of $\Delta-c$ with the same eigenvalue $z$ and $T_t\ell_y(\psi(z))(x)=e^{-zt}\ell_y(\psi(z))(x)$. Since $\Re z>0$, it follows that $T_t\ell_y(\psi(z))(x)\to 0$ as $t\to \infty$. Therefore $V_1\subset B_0$. Since $B_0$ is a vector space we have  $\mathrm{span}(V_1)\subset B_0$.
We recall that $(C_0(S))^\ast$ is the set of complex measures.
If a nonzero measure  $\mu\in (C_0(S))^\ast$ annihilates $\mathrm{span}(V_1)$, then for every $z\in \Omega_\infty^+$, $\mu\ast \phi_{\Gamma(z)}\equiv 0$. In other words for every  $\lambda\in \Gamma(\Omega_\infty^+)$, $\mu\ast \phi_\lambda\equiv 0$. Noting that for any fixed $y\in S$, $\lambda\mapsto \mu\ast\phi_\lambda(y)$ is an analytic function on $S_1^\circ$ and that $\Gamma(\Omega_\infty^+)$ is an open set we conclude that $\mu\ast \phi_\lambda(y)=0$ for all $\lambda\in S_1$ and for all $y\in S$. This implies that $\mu=0$.
  Hence $\mathrm{span}\, V_1$ as well as $B_0$ is dense in $C_0(S)$.

We define $V_2=\{\ell_y(\psi(z))(x)\mid z\in \Omega_\infty^-, y\in S \}$ and $V_3=\{\ell_y(\psi(z))(x)\mid z\in \Omega_\infty\cap i\Q,  y\in S \}$. Again the elements $\ell_y(\psi(z))(x)$ of $V_2$ and $V_3$ are eigenfunctions of $\Delta-c$ with eigenvalue $z$. Argument analogous to (i) now shows that $\mathrm{span}\, V_2\subset B_\infty$ and $\mathrm{span}\, V_3\subset B_{\mathrm{Per}}$.  Finally  the  denseness of  $B_\infty$ and $B_{\mathrm{Per}}$ in $C_0(S)$ will follow through  the argument  used above to show denseness of  $\mathrm{span}\, V_1$.   We omit the details to avoid repetition.

Argument for  the converse  statement is also same as that of (i). (See Remark \ref{non-subspace-chaotic}.)
\vspace{.1in}

(iii)  First we deal with the case $2<p<\infty$ and $1\le q<\infty$.
To show that $T_t$  is not hypercyclic we take a nonzero  $\phi\in L^{p',q'}(S)$.  We claim that $\{T_t^\ast \phi\mid t\ge 0\}$
is a bounded set in $L^{p',q'}(S)$. We note that $T_t^\ast \phi=T_t\phi= \phi\ast p_t$ where $p_t=e^{ct}h_t$. From Lemma~\ref{herz}
we see that $\|T_t^\ast\phi\|_{p', q'}\le \what{p_t}(i\gamma_p\rho)\|\phi\|_{p',q'}=e^{(c-c_p)t} \|\phi\|_{p',q'}$   As $c\le c_p$ the claim is established.
Therefore by  Proposition~\ref{Not-chaotic-hyper}~(iii), $T_t$ is not hypercylic on $L^{p,q}(S)$.

For the case  $2<p<\infty$ and $q=\infty$,  we take $\phi\in L^{p', 1}(S)$. Then $\phi\in (L^{p, \infty})^\ast(S)$. By Lemma~\ref{herz},  $\|T_t^\ast\phi\|_{p', 1}=\|T_t\phi\|_{p', 1} \le e^{(c-c_p)t} \|\phi\|_{p',1}$. Then \[\|T_t\phi\|_{(L^{p, \infty}(S))^\ast}=\sup_{\psi\in L^{p, \infty}(S)}\frac{|\langle T_t\phi, \psi\rangle|}{\|\psi\|_{p, \infty}}\le \frac{e^{(c-c_p)t} \|\phi\|_{p',1}\|\psi\|_{p, \infty}}{\|\psi\|_{p, \infty}}=e^{(c-c_p)t} \|\phi\|_{p',1}.\] Rest of the argument is same as the previous case. The case of $L^\infty(S)$ can be treated analogously,  taking $\phi\in L^1(S)\subset (L^\infty(S))^\ast$. From Remark \ref{non-subspace-chaotic} it is clear that $T_t$ cannot be subspace-chaotic on these spaces.
\end{proof}

\begin{proof}[proof of Theorem~\ref{result-S-2}]
(i) Once we notice that for $1\le q\le \infty$, the $L^{2, q}$-spectrum of $\Delta-c$ is either empty or lies entirely on $\R$, it follows  from   Proposition \ref{Not-chaotic-hyper} (i) that $T_t$ is not chaotic or subspace-chaotic.

For the last part of (i), the argument is similar to  what is used for the proof of  Theorem~\ref{result-S-1}~(iii). We can show that  for a nonzero function $\phi\in L^{2, q'}(S)$, $\|T_t^\ast \phi\|_{2, q'}\le e^{(c-\rho^2)t}\|\phi\|_{2, q'}$. Similarly for a nonzero function
$\phi\in L^{2, 1}(S)$, $\|T_t^\ast \phi\|_{(L^{2, \infty}(S))^\ast}\le e^{(c-\rho^2)t}\|\phi\|_{2, 1}$. Proposition~\ref{Not-chaotic-hyper}~(iii) now proves the assertion as $c\le \rho^2$.

(ii)
To establish the non-hypercyclicity,  (in view of Proposition~\ref{Not-chaotic-hyper}~(ii)) it suffices to show that the  point spectrum of $(\Delta-c)^\ast=(\Delta-c)$ is nonempty on the dual spaces of these spaces, which is indeed the case as: $\phi_\lambda\in L^{2, \infty}(S)$ for any nonzero  $\lambda\in \R$; $\phi_\lambda\in (L^{p,q}(S))^\ast=L^{p', q'}(S)$ for $\lambda\in S_p^\circ$ (see section 3) and  $\phi_\lambda\in L^\infty(S)$ if $\lambda\in S_1$. This also shows that $T_t$ is not subspace-chaotic on these spaces.
\end{proof}

\section{Existence of periodic points of $T_t$}
As in the previous sections let $T_t=e^{-t(\Delta-c)}$ for $t\ge 0$, where  $c\in \R$ is fixed and for any $p\ge 1$, $c_p=4\rho^2/pp'$. It is clear that if $(\Delta-c)f=0$ for some suitable function $f$, then $T_tf=f$ for all $t\ge 0$. However from the hypothesis $T_tf=f$ for some $t>0$, it does not seem to be straightforward to conclude that $(\Delta-c)f=0$ (see \cite{DSW}). This forces us to go through a round-about argument involving Wiener Tauberian theorem (WTT) to prove Proposition~\ref{no-periodic} stated below, which is the aim of this section.
We begin with the following observations.
\begin{enumerate}
\item[(a)] For $2<p\le \infty$ and $c\ge c_p$, $T_t$ has periodic  points in $L^{p, \infty}(S)$.
\item[(b)] For $c>\rho^2$, $T_t$ has periodic point on $L^{2, \infty}(S)$.
\item[(c)] For any $c\in \R$, $T_t$ has no periodic point in the spaces $L^1(S)$,  $L^{2,r}(S)$ with $1\le r\le 2$ and  $L^{p, q}(S)$ with $1<p<2, 1\le q\le \infty$.
\end{enumerate}
Indeed,  the assertion (a) with  $c>c_p$ is proved in Theorem~\ref{result-S-1}~(i),~(ii) since $L^{p,q}(S)\subset L^{p, \infty}(S)$ and $L^{\infty, \infty}(S)=L^\infty(S)$. For the case $c=c_p$ in (a),  we recall that $\phi_{i\gamma_p\rho}\in L^{p, \infty}(S)$.  It can be verified easily that for all $t\ge 0$, $T_t \phi_{i\gamma_p\rho}=\phi_{i\gamma_p\rho}$. In (b) if $\lambda\in \R$ satisfies $\lambda^2+\rho^2=c$ then $\lambda\neq 0$ and hence  $\phi_\lambda\in L^{2, \infty}(S)$ (see section 3). It can be verified that $T_t\phi_\lambda=\phi_\lambda$.
For (c) we first note that $L^{2,r}(S)\subset L^2$ for $1\le r\le 2$. Therefore it is enough to consider the spaces $L^1(S), L^2(S)$ and $L^{p,q}(S)$ with $p,q$ in the given range. We recall that if $f$ is in one of these spaces then its Fourier transform $\wtilde{f}(\lambda, n)$ exists as a measurable function in $\lambda\in \R$ for almost every $n\in N$. Therefore if for some $t>0$, $T_tf-f=0$ where  $f$ is in any of these spaces, then taking Fourier transform we have $(e^{-t(\lambda^2+\rho^2-c)}-1)\wtilde{f}(\lambda, n)=0$ for almost every $\lambda\in \R$ and almost every $n\in N$. From this and (\ref{inversion}) it follows that $f=0$.

\subsection{WTT and its application}
Some versions of the WTT for radial (i.e. $K$-biinvariant) functions in the Lorentz spaces $L^{p,q}(G/K)$, $1<p<2, 1\le q<\infty$  and in $L^1(G/K)$ were established (see \cite[Theorem 6.1, Remark 6.1.1]{PRS}) for the rank one symmetric spaces. It is not difficult to see that exactly the same argument (which uses the result of disk algebra following \cite[Theorem 1.1]{BW}), yields the corresponding results for the radial functions on DR spaces, which we shall state below.

We recall that (see \cite{KRS}) for $1\le p<2$, $L^{p,1}(S)\ast L^{p,1}(S)^\#\subset L^{p,1}(S)$ and $\|f\ast g\|_{p,1}\le C\|f\|_{p,1}\|g\|_{p,1}$ where $f\in L^{p,1}(S)$ and $g\in L^{p,1}(S)^\#$. In particular $L^{p,1}(S)^\#$ is a Banach algebra under convolution.
\begin{proposition} \label{WTT} Let $1<p<2$ be fixed. Suppose that for a radial function  $f\in L^{p,1}(S)$ {\em(}respectively  $f\in L^p(S), 1\le p<2${\em)}, its spherical Fourier transform $\what{f}$ satisfies,
\begin{enumerate}
\item[(i)] $\what{f}$ extends analytically on $S_p^\varepsilon=\{z\in \C\mid |\Im z|< |\gamma_p|+\varepsilon\}$ for some $\varepsilon>0$;

\item[(ii)] $\lim_{|\lambda|\to \infty} \what{f}(\lambda)=0$ on $S_p^\varepsilon$;

\item[(iii)] for all $\lambda\in S_p^\varepsilon$, $\what{f}(\lambda)\neq 0$;

\item[(iv)] for  $t\in \R$, $\limsup_{|t|\rightarrow\infty}|\what{f}(t)e^{Ke^{|t|}}|>0$  for all $K>0$.
\end{enumerate}
Then  the ideal {\em(}respectively the  $L^1(S)^\#$-module{\em)} generated by $f$ is dense in
$L^{p,1}(S)^\#$ {\em(}respectively in $L^p(S)^\#${\em)}.
\end{proposition}
Following extension is immediate: if  a radial function $f\in L^{p,1}(S), 1\le p<2$ satisfies conditions (i)-(iv), then the left $L^{p,1}(S)$-module generated by $f$  is dense in $L^{p,1}(S)$. Let us denote by $M$ the closed left $L^{p,1}(S)$-module generated by $f$ in $L^{p,1}(S)$. It is clear from Proposition \ref{WTT} that
$M\supset L^{p,1}(S)^\#$, hence in particular $h_t\in M$ for all $t\ge 0$. We take any  $g\in L^{p,1}(S)$. Then $g\ast h_t\in M$. Since  $g\ast h_t\to g$ in $L^{p,1}(S)$ and $M$ is closed, we have $g\in M$. The  same  argument also shows that if a  radial function $f\in L^{p}(S), 1\le p<2$  satisfies conditions (i)-(iv), then the left $L^{1}(S)$-module generated by $f$ in $L^{p}(S)$ is dense in $L^{p}(S)$.
We shall now apply these results to prove the following.
\begin{proposition}\label{no-periodic}
For $p\in (2, \infty)$,  $q\in [1, \infty]$ and $c<c_p$, $T_t$ has no periodic point in $L^{p,q}(S)$. If $c<c_\infty=0$, then $T_t$ has no periodic  point in $L^\infty(S)$.
\end{proposition}
\begin{proof}
Since $L^{p, q}(S)\subset L^{p, \infty}(S)$ for all $q$, it suffices to show  that the assumption $T_tf=f$ for  a nonzero  $f\in L^{p, \infty}(S)$ and a $t>0$, leads to a contradiction. Indeed, this assumption implies that $e^{ct} f\ast h_t - f=0$. Convolving with $h_{t'}$ for some $t'> 0$, we get
$f\ast (e^{ct}h_{t+t'}-h_{t'})=0$. Let $h=(e^{ct}h_{t+t'}-h_{t'})$.  Then $h\in L^{p',1}(S)$, $h$ is radial and $f\ast h=0$.  For $0<\varepsilon<c$,    $\what{h}(\lambda)=e^{-t'(\lambda^2+\rho^2)}(e^{-t(\lambda^2+\rho^2-c)}-1)$ extends analytically to $S_p^\varepsilon$ and $\what{h}(\lambda)\neq 0$ for all $\lambda\in S_p^\varepsilon$. Indeed, $\what{h}$ satisfies conditions (i)-(iv) of Proposition~\ref{WTT}.  It follows from Proposition~\ref{WTT} and the subsequent discussion,  that the left $L^{p', 1}$-module  generated by $e^{ct}h_{t+t'}-h_{t'}$  is  dense in  $L^{p',1}(S)$.  Since   for any $g_1\in L^{p',1}(S)$, $\langle f, g_1\ast h\rangle=\langle f\ast h, g_1\rangle=0$,  we have for any $g\in L^{p',1}(S)$, $\langle f, g\rangle=0$. This implies that $f=0$, since $f\in L^{p, \infty}(S)=(L^{p',1}(S))^\ast$.
For the case of $L^\infty(S)$,  we argue the same way and use Proposition~\ref{WTT} for $L^1(S)$.
\end{proof}
\vspace{.1in}

\noindent{\bf Acknowledgements:} The author would like to thank S. K. Ray for making him aware of  the  work of Ji and Web and for numerous suggestions.
He is also thankful to S. C. Bagchi for many conversations.

  \end{document}